\documentclass[reqno,12 pt]{amsart}
\usepackage{amssymb}

\def \1{\bf 1}

\def \F{\mathbb F}
\def \Z{\mathbb Z}
\def \C{\mathbb C}
\def \R{\mathbb R}

\def \l{\lambda}

\def \wt{{\rm wt}}

\def \mod{{\rm mod}}

\def \<{\langle}
\def \>{\rangle}

\def \o{\omega}

\def \ch{{\rm ch}}
\def \a{\alpha }

\def \pf{\noindent {\bf Proof: \,}}

\def \1{{\bf 1}}
\def\Ve{V^{0}}
\def\ha{\frac{1}{2}}
\def\se{\frac{1}{16}}

\newcommand{\abs}[1]{\left\vert#1\right\vert}

\newcommand{\lan}{\langle}
\newcommand{\ran}{\rangle}

\newcommand{\be}{\beta}
\newcommand{\al}{\alpha}

\newcommand{\supp}{\mathrm{supp}}
\newcommand{\sC}{\mathcal{C}}
\newcommand{\sD}{\mathcal{D}}
\newcommand{\sH}{\mathcal{H}}

\newtheorem{thmm}{Theorem}
\newtheorem{thm}{Theorem}[section]
\newtheorem{prop}[thm]{Proposition}
\newtheorem{cor}[thm]{Corollary}
\newtheorem{lem}[thm]{Lemma}
\newtheorem{rem}[thm]{Remark}
\newtheorem{defn}[thm]{Definition}

\begin{document}
\title{On the uniqueness of the moonshine vertex operator algebra}
%\bsays{I did not want to leave ``FLM'' in the title since that may be useless for searches and is too informal.  If I put full names at the end, the running head will be too long.  The full names are given several times in the text.    }  

\author{Chongying Dong}
\address{Department of Mathematics, University of California,
Santa Cruz, CA 95064}
\thanks{C.D. is supported by NSF grants, China NSF grant 10328102 and faculty research funds granted by
the University of California at Santa Cruz}

\author{Robert L. Griess Jr.}
\address{Department of Mathematics, University of
Michigan, Ann Arbor, MI 48109-1109}
\thanks{R.L.G. is supported by NSA grant USDOD-MDA904-03-1-0098}
\author{Ching Hung Lam}
\address{Department of Mathematics, National Cheng Kung University, Tainan, Taiwan 701}
\thanks{C.L. is supported by NSC grant 93-2115-M-006-012 of Taiwan
and National Center for Theoretical Sciences, Taiwan}

\begin{abstract} It is proved that the vertex operator algebra  $V$
is isomorphic to the moonshine VOA  $V^{\natural}$ of
Frenkel-Lepowsky-Meurman if it satisfies conditions (a,b,c,d) or
(a$^\prime$,b,c,d).  These conditions are:

(a)  $V$ is the only irreducible module for itself and $V$ is
$C_2$-cofinite;

(a$^\prime$) $\dim V_n\leq \dim V^{\natural}_n$ for $n\geq 3$;

(b) the central charge is 24;

(c) $V_1=0$;

(d) $V_2$ (under the first product on $V$) is isomorphic to the
Griess algebra.

\noindent Our two main theorems therefore establish a  weak
version of the FLM uniqueness conjecture for the moonshine vertex
operator algebra.  We believe that these are the first such
results.
\end{abstract}

\maketitle

\section{Introduction}

The moonshine vertex operator algebra $V^{\natural}$ constructed
by Frenkel-Lepowsky-Meurman \cite{FLM1},\cite{FLM2} not only
proves a conjecture by McKay-Thompson but also  plays a
fundamental role in shaping the theory of vertex operator algebra.
In the introduction of \cite{FLM2}, Frenkel-Lepowsky-Meurman
conjectured that the $V^{\natural}$ can be characterized by the
following three conditions:

(a) the VOA $V^{\natural}$ is the only irreducible ordinary module
for itself;

(b) the central charge of $V^{\natural}$ is 24;

(c) $V^{\natural}_1=0.$

We call their conjecture {\it the Frenkel-Lepowsky-Meurman conjecture. }
These conditions are natural analogues of conditions which
characterize the binary Golay code and the Leech lattice.

Conditions (b) and (c) are clear from the construction. Condition
(a) is proved in \cite{D} by using the 48 commuting Virasoro
elements of central charge $\frac 1 2$ discovered in \cite{DMZ}.
Furthermore, $V^{\natural}$ is rational \cite{DLM2},\cite{DGH}.
Although the theory of vertex operator algebra has developed a lot
since \cite{FLM2}, including some uniqueness results for certain
VOAs \cite{LX}, \cite{DM2}, \cite{DM3},  there has been  no real
progress in proving their  conjecture.

In this paper we prove two  weak versions of the Frenkel-Lepowsky-Meurman  conjecture:
\begin{thmm}\label{mt1} Let $V$ be a $C_2$-cofinite vertex operator algebra satisfying
(a)-(c). We also assume that $V_2$ is isomorphic to the Griess
algebra.
 Then $V$ is isomorphic to $V^{\natural}.$
\end{thmm}

In the second main theorem, we replace condition (a)  by the
assumption that $\dim V_n\leq \dim V^{\natural}_n$ for $n\geq 3.$
\begin{thmm}\label{mt2} Let $V$ be a simple vertex operator algebra satisfying
(b)-(c). We also assume that $V_2$ is isomorphic to the Griess
algebra and $\dim V_n\leq \dim V^{\natural}_n$ for $n\geq 3.$
 Then $V$ is isomorphic to $V^{\natural}.$
\end{thmm}

We now discuss the theorems and background. The weight two
subspace $V^{\natural}_2$ of $V^{\natural}$ with the product which takes the pair $u, v$ to 
$u_1v$, where $u_1$ is the component operator of the
vertex operator $Y(u,z)=\sum_{n\in\Z}u_nz^{-n-1}$ \cite{FLM2}, is
the Griess algebra \cite{G}, which is a commutative nonassociative
algebra of dimension 196884. Moreover, $V^{\natural}$ is generated
by $V^{\natural}_2$ and $V^{\natural}$ is an irreducible module
for the affinization of the Griess algebra \cite{FLM2}. So, in
order to understand the moonshine vertex operator algebra, one
must know the Griess algebra and its affinization very well. It
seems that a complete proof of FLM's uniqueness conjecture needs a
better understanding of the Griess algebra. Unfortunately, there
does not yet exist a characterization of the Griess algebra
(independent of its connection to the monster simple group).
Also, the affinization of the Griess algebra is not a Lie algebra
and lacks a highest weight module theory. From this point of view,
$V^{\natural}$ is a very difficult vertex operator algebra.

The study of the moonshine vertex operator algebras in terms of
minimal series of the Virasoro algebras was initiated in
\cite{DMZ}. This is equivalent to the study the maximal
associative subalgebra of the Griess algebra. In \cite{DMZ}, we
find 48 mutually commutative Virasoro algebras with central charge
$\frac 1 2.$ As a result, a tensor product $T_{48}$ of 48 vertex
operator algebras, associated to the highest weight unitary
representations of the Virasoro algebra with central charge $\frac
1 2$ is a subalgebra of $V^{\natural}$ and $V^{\natural}$
decomposes into a direct sum of finitely many irreducible modules
for $T_{48}$ as $T_{48}$ is rational and the homogeneous summands for $V$ are finite dimensional. A lot of progress on
the study of the moonshine vertex operator algebra has been made
by using the subalgebra $T_{48}$ and vertex operator subalgebras
associated to the other minimal unitary series for the Virasoro
algebras \cite{DLMN}, \cite{DGH}, \cite{KLY}, \cite{M3}. The
discovery of the $T_{48}$
inside $V^{\natural}$ also inspired
the study of code vertex operator algebras and framed vertex
operator algebras \cite{M2}, \cite{DGH}.

A {\it frame} in $V^{\natural}$ is a set of 48 mutually orthogonal Virasoro
elements with central charge $\frac 1 2.$ The subalgebra $T_{48}$
depends on a frame as studied in \cite{DGH}. It is proved in
\cite{DGH} that for any choice of 48 commuting Virasoro algebras
there are two codes $C$ and $D$ associated to the decomposition of
$V^{\natural}$ into irreducible $T_{48}$-modules. Each irreducible
$T_{48}$-module is a tensor product of 48 unitary highest weight
modules   $L(\frac{1}{2},h)$ for the Virasoro algebra with central
charge $\frac 1 2$ where $h$ can take only three values
$0,\frac{1}{2}, \frac{1}{16}.$ The code $C$ tells us the
irreducible $T_{48}$-modules occurring in $V^{\natural}$ which are
a tensor product of $L(\frac{1}{2},h)$ for $h=0$ or $\frac{1}{2}.$
Similarly, the code $D$ indicates the appearance of irreducible
$T_{48}$-modules whose tensor factors have at least one
$L(\frac{1}{2},\frac{1}{16}).$   The fusion rules for the vertex
operator algebra $L(\frac{1}{2},0)$  indicate that we should
consider a frame so that $C$ has maximal possible dimension and
$D$ minimal possible dimension. These respective dimensions are 41
and 7. The reason for using our particular frame is that, for the
code VOA which arises, {\it all irreducible modules are simple
currents} (see Theorem \ref{fusion} in this paper). The uniqueness
of $V^\natural$ then  follows from known uniqueness results for
certain smaller VOAs, those which are simple current extensions of
code VOAs.

The main strategy in proving the theorem is to use this particular
frame. Since we assume that the weight 2 subspace of the abstract
vertex operator algebra in the theorem is isomorphic to the Griess
algebra, we can use the theory of framed vertex operator algebra
developed in \cite{DGH} and \cite{M2} to investigate the structure
of such vertex operator algebras.

Although we assume that $V_2\cong V_2^\natural$ (as algebras), we
can not claim automatically that any VF in $V^\natural$
corresponds to a VF in $V$.   The difficult point is to prove that
a Virasoro vector in $V_2^\natural$ generates a subVOA which is
{\it simple}, i.e. an irreducible highest weight module. This is
where we make use of the other assumptions in our main theorems.
The proof involves both character theory for the Virasoro algebra
with central charge $\frac 1 2$ and an explicit expression for the
$J$-function.

It seems that there is still a long way to go to settle the  FLM
conjecture. The main difficulty is that we do not have much theory
of finite dimensional commutative nonassociative algebras which
could be applicable to a 196884-dimensional degree 2 summand of a
VOA satisfying our conditions (a,b,c) (see \cite{G1}).   In a
sense, this paper reduces the uniqueness of the moonshine vertex
operator algebra to the uniqueness of the Griess algebra.

\section{Notations}

Most of our notations are fairly standard in the VOA literature.
For the reader, we note a few below.

\bigskip

codes $C=C(F), D=D(F)$: see Section 4;

codes $\sC, \sD$: see Section 6;

$j(q), J(q)$ :  the elliptic modular function and the elliptic
modular function with constant term set equal to 0, i.e.,
$J(q)=j(q)-744$;

$\lan \omega_i \ran$    :   the subVOA generated by $\omega_i$;

VF  :   Virasoro frame, see Section 4;

$Vir(\omega_i)$  :  the Virasoro algebra spanned by the modes of
the Virasoro element $\omega_i$ and the scalars;

$V^I$ : see Section 4;

$V^0$ or $V^{\emptyset}$ : the case of $V^I$ for $I=0$ or
$\emptyset$;

$V^\natural$ :  the moonshine VOA, constructed in \cite{FLM2};

$(V^\natural)^0$  :   this is $V^0$ for $V=V^\natural$.

\section{Various modules for vertex operator algebras}
\setcounter{equation}{0}

Let $(V,Y,\1,\omega)$ be a vertex operator algebra. We recall
various notion of modules (cf. \cite{FLM2}, \cite{DLM1}).

\bigskip

A {\em weak} $V$ module is a vector space $M$ with a linear map
$Y_M:V \rightarrow End(M)[[z,z^{-1}]]$ where $v \mapsto
Y_M(v,z)=\sum_{n \in \Z}v_n z^{-n-1}$, $v_n \in End(M)$. In
addition $Y_M$ satisfies the following:

1) $v_nw=0$ for $n>>0$ where $v \in V$ and $w \in M$;

2) $Y_M( {\textbf 1},z)=Id_M$;

3) The Jacobi Identity
\begin{eqnarray*}
& &z^{-1}_0\delta\left(\frac{z_1-z_2}{z_0}\right)
Y_M(u,z_1)Y_M(v,z_2)-z^{-1}_0\delta\left(\frac{z_2-z_1}{-z_0}\right)
Y_M(v,z_2)Y_M(u,z_1)\\
& & \ \ \ =z_2^{-1}\delta\left(\frac{z_1-z_0}{z_2}\right)
Y_M(Y(u,z_0)v,z_2)
\end{eqnarray*}
holds.

\bigskip

An {\em admissible} $V$ module is a weak $V$ module which carries
a $\Z_+$-grading, $M=\bigoplus_{n \in \Z_+} M(n)$, such that
 $v_m M(n) \subseteq M(n+\wt v-m-1)$

\bigskip

An {\em ordinary} $V$ module is a weak $V$ module which carries a
$\C$-grading, $M=\bigoplus_{\l \in \C} M_{\l}$, such that:

1) $dim(M_{\l})< \infty$ for all $\l \in \C$;

2) $M_{\l+n}=0$ for fixed $\l$ and $n<<0$ (depending on $\l$);

3) $L(0)w=\l w={\rm wt}(w) w$, for $w \in M_\l$.
\bigskip

It is easy to prove that an ordinary module is admissible.

A vertex operator algebra is called {\em rational} if every
admissible module is a direct sum of simple admissible modules.
That is, a VOA is rational if there is complete reducibility of
the category of admissible modules. It is proved in \cite{DLM2}
that if $V$ is rational there are only finitely many irreducible
admissible modules up to isomorphism and each irreducible
admissible module is ordinary.
% \bsays{Is there a Krull Schmidt
%theorem for rational VOAs.  It must be so if $V$ has a highest
%weight theory, but what about in general?  }  \chosays{What is
%Krull Schmidt theorem?}
%\bsays{ It is more correct to say Remak-Krull-Schmidt.  Let S be a set and Mod a category of S-modules.  An object of Mod is a group (not necessarily abelian) with a set map from S to the set of endomorphisms of that group.  Subobjects and quotients must involve S-invariant subgroups.  The classic RKS Theorem says that if M is an object of Mod and M satisfies ACC and DCC on S-subgroups then there exists a direct sum M=M1 +...+Mp into S-indecomposables.  Also if N1+...+Nq =M is another, then p=q and up to reindexing Mi and Ni are isomorphic.
%\smallskip
%There is a RKS theorem for f.g. modules over a PID even though these modules fail to satisfy DCC.
%\smallskip
%There is  a RKS theorem for f.g. highest weight modules for a ss Lie algebra.  I think also for a KMLA.
%\smallskip
%There is not a RKS theorem for f.g modules for integer group rings of finite groups.  I think Quat(16) gives a counterexample.
%\smallskip
%There is not a RKS theorem for f.g. free modules over noncommutative rings in general.  Consider End(V) where V has countably infinite dimension and take a linear isomorphism of V with V+V.
%\smallskip
%Hence my question: is there a RKS theorem for f.g. modules over a rational VOA.  I guess probably yes if there is a highest weight theory, but no guess in general.  }  \chosays{I dont think it is related to this paper.}

A vertex operator algebra $V$ is called {\em holomorphic} if it is
rational and the only irreducible ordinary module is itself. In
this case $V$ is also the only irreducible admissible module.

A vertex operator algebra is called {\em regular} if every weak
module is a direct sum of simple ordinary modules. So, regularity
implies rationality.

A vertex operator algebra $V$ is called {\it $C_2$-cofinite} if
$V/C_2(V)$ is finite dimensional where $C_2(V)=\<u_{-2}v|u,v\in
V\>.$

\section{Framed vertex operator algebras}
\setcounter{equation}{0}

In this section we review the framed vertex operator algebras and
related results from \cite{DMZ} and \cite{DGH}.

Let $L(c,h)$ be the irreducible highest weight module for the
Virasoro algebra with central charge $c$ and highest weight $h.$
The $L(\frac 1 2 ,0)$-module  $L(\frac{1}{2},h)$ is unitary if and
only if $h=0,\frac{1}{2},\frac{1}{16}$ \cite{FQS}, \cite{GKO}.
Moreover, $L(\frac{1}{2},0)$ is a rational vertex operator algebra
and $L(\frac{1}{2},h)$ for $h=  0,\frac{1}{2},\frac{1}{16}$ gives
a complete list of inequivalent irreducible
$L(\frac{1}{2},h)$-modules.

We first recall the notion of framed vertex operator algebra. Let
$r$ be a nonnegative integer. A {\it framed vertex operator
algebra} (FVOA)  is a simple vertex operator
algebra $(V,Y,\1,\o)$  
satisfying the following conditions:  there exist
$\omega_i\in V$ for $i=1$,~$\ldots$,~$r$ such that  (a) each
$\o_i$ generates a copy of the simple Virasoro vertex operator
algebra $L(\frac{1}{2},0)$ of central charge $\frac{1}{2}$ and the
component operators $L^i(n)$ of
$Y(\omega_i,z)=\sum_{n\in\Z}L^i(n)z^{-n-2}$ satisfy
$[L^i(m),L^i(n)]=(m-n)L^i(m+n)+\frac{m^3-m}{24}\delta_{m,-n};$ (b)
The $r$ Virasoro algebras $Vir(\omega_i )$, spanned by the modes
of $Y(\omega_i ,z)$ and the identity, are mutually commutative;
and (c) $\omega=\omega_1+\cdots+\omega_{r}$. The set
$\{\o_1,\ldots,\o_r\}$ is called a  {\it Virasoro frame}  (VF).

>From now on we assume that $V$ is a FVOA of central charge
$\frac{r}{2}$  with frame $F:=\{\o_1,\ldots,\o_r\}$.   Let $T_r$
be the vertex operator algebra generated by $\omega_i$ for
$i=1,...,r.$ Then $T_r$ is isomorphic to
$L(\frac{1}{2},0)^{\otimes r}$ and  its irreducible modules  are the 
$L(h_1,...,h_r):=L(\frac{1}{2},h_1)\otimes\cdots \otimes
L(\frac{1}{2},h_r)$ for $h_i=0,\frac{1}{2},\frac{1}{16}.$ Since
$T_r$ is a rational vertex operator algebra, $V$ is a completely
reducible $T_r$-module. That is,
\begin{equation}\label{2.1}
V \cong \bigoplus_{h_i\in\{0,\frac{1}{2},\frac{1}{16}\}}
m_{h_1,\ldots, h_{r}}L(h_1,\ldots,h_{r})
\end{equation}
where the nonnegative integer $m_{h_1,\ldots,h_r}$ is the
multiplicity of $L(h_1,\ldots,h_r)$ in $V$. In particular, all the
multiplicities are finite and $m_{h_1,\ldots,h_r}$ is at most $1$
if all $h_i$ are different from $\frac{1}{16}$.

There are two binary codes $C=C(F)$ and $D=D(F)$ associated to the
decomposition (\ref{2.1}). In order to define the code $D$ we
identify a subset $I$ of $\{1,...,r\}$ with a codeword
$d=(d_1,...,d_r)\in \F_2^r$ where $d_i=1$ if $i\in I$ and $d_0=0$
elsewhere. Let $I$ be a subset of $\{1,\ldots,r\}$. Define $V^I$
as the sum of all irreducible submodules isomorphic to one of the
irreducibles $L(h_1,\ldots , h_r)$ such that $h_i=\frac{1}{16}$ if
and only if $i \in I$. Then
$$V=\bigoplus_{I\subseteq \{1,\ldots,r\}}V^I.$$
Set
\begin{equation}\label{2.2}
D=D(F):=\{I\in \F_2^r \mid V^I \ne 0\}.
\end{equation}

For $c=(c_1,...,c_r)\in \F_2^r$, we define
$V(c)=m_{h_1,...,h_r}L(h_1,...,h_r)$ where $h_i=\frac{1}{2}$ if
$c_i=1$ and $h_i=0$ elsewhere. Set
 \begin{equation}\label{2.3}
C=C(F):=\{c\in \F_2^r \mid V(c) \ne 0\}.
\end{equation}
Then $V^{\emptyset}=V^0=\bigoplus_{c\in C}V(c)$.

Here we summarize the main result about FVOAs   from \cite{DGH}
\begin{thm}\label{DGH} Let $V$ be a FVOA. Then

(a) $V=\oplus_{n\geq 0}V_n$ with $V_0=\C \1.$

(b) $V$ is rational.

(c) $C$ and $D$ are binary codes and
$$C\subset  D^{\perp}=\{x=(x_1,...,x_r)\in \F_2^r|x\cdot d=0
\forall d\in D\}.$$ Moreover, $V$ is holomorphic if and only if
$C=D^{\perp}.$

(d) $\Ve$ is a simple vertex operator algebra and the $V^I$ are
irreducible $\Ve$-modules. Moreover $V^I$ and $V^J$ are
inequivalent if $I\ne J$.

(e) For any $I,J\in D$ and $0\ne v\in V^J$ we have
$V^{I+J}=span\{u_nv|u\in V^I,n\in\Z\}.$

(f) Let $I \subseteq \{1,\ldots,r\}$ be given and suppose that
$(h_1,\ldots,h_r)$ and $(h_1',\ldots,h'_r)$ are $r$-tuples with
$h_i$, $h_i'\in\{0,\ha,\se\}$ such that $h_i=\frac{1}{16}$
(resp.~$h_i'=\frac{1}{16}$) if and only if $i\in I$. If both
$m_{h_1,\ldots,h_r}$ and $m_{h_1',\ldots,h'_r}$  are nonzero then
$m_{h_1,\ldots,h_r}=m_{h_1',\ldots,h'_r}$. That is, all
irreducible modules inside $V^I$ for $T_r$ have the same
multiplicities.

(g) For any $c,d\in C$ and  $0\ne v\in V(d)$ we have
$V(c+d)=span\{u_nv|u\in V(c),n\in\Z\}.$
\end{thm}

\section{Code VOA $M_C$}

In this section we review and extend results on code VOAs and
their modules, following \cite{M1}-\cite{M3} and \cite{La}.

We shall sometimes consider an integer modulo 2 as its Euclidean
lift, i.e., its representative 0 or 1 in $\Z$, so that  when $\al
\in \Z_2$, $\frac 1 2 \al$ makes sense as the rational number 0 or
$\frac 1 2$.

Let $C$ be an even binary code. For any $\al= (\al_1, \dots,
\al_n)\in C$, denote
\[
M_\al= L(\frac{1}2,\frac{\al_1}2) \otimes \cdots \otimes
L(\frac{1}2,\frac{\al_n}2) \quad \text{ and }\quad M_C=
\bigoplus_{\al\in C} M_\al.
\]
Note that $M_C$ is a simple current extension of
$T_n=L(\frac{1}2,0)^{\otimes {n}}$ and it has a unique VOA
structure over $\C$ (cf. \cite{DM2}, \cite{M2}).  This will be
used to deduce the uniqueness of $V^\natural$.

\begin{rem} We use $M_C$ for a code VOA instead of $M_D$ given
in \cite{M2} in this paper. This is consistent with our code $C$
defined in Section 3. In fact, $M_C$ is a framed VOA with frame
$F$ satisfying $C(F)=C$ and $D(F)=0.$
\end{rem}

\begin{rem}
For any $\beta\in \mathbb{Z}_2^n$, one can define an automorphism
$\sigma_{\beta}:M_C \to M_C$ by
\[
\sigma_{\beta}(u)=(-1)^{\lan \alpha, \beta \ran}u\qquad \text{ for
} u\in M_{\alpha}.
\]
 This automorphism is called a coordinate  automorphism.
Note that $\sigma_{\beta}=\sigma_{\beta'}$ if and only if
$\beta+\beta'\in C^{\bot}$ and the subgroup $P$ generated by $\{
\sigma_{\beta}|\ \beta \in \Z_2^n\}$ is isomorphic to $\Z_2^n/
C^\perp$. Moreover, the fixed subalgebra ${M_C}^P$ is $T_n$ (cf.
\cite{M1}).
\end{rem}

We first study the representations of the code VOA $M_C.$ Let $W$
be an irreducible $M_{C}$-module. Then $W$ can be written as a
direct sum of irreducible $T:=T_n$-modules,
\[
W\cong \bigoplus_{h_{i}\in \{0,\frac{1}2,\frac{1}{16}\}} m_{h_{1},
\dots,h_{n}} L(h_1,\cdots, h_n).
\]

\begin{defn}
Define $\tau(L(h_1,\cdots, h_n))=(a_{1},\cdots ,%
a_{n})$ $\in \mathbb{Z}_{2}^{n}$ such that
\[
a_i=\left\{
\begin{array}{lll}
  0 & \text{\quad if }& h_{i}=0\text{ or }\frac{1}{2} \\
  1 & \text{\quad if }&h_{i}=\frac{1}{16}
\end{array}
\right. .
\]
This binary word is called the $\tau $-word of $L(h_1,\cdots, h_n).$
%\bsays{I don't understand why we use $\tilde\tau$ and not just
%${\tau}$. } \chisays{ `` $\tau$ " is usually used to denote
%Miyamoto's involution. So, we often use $\tilde{\tau}$ to denote the
%$\tau$-word. I think we can use $\tau$ here also.}\bsays{I see
%little danger of confusion, so may we use only $\tau$?}  \chosays{I
%think $\tilde{\tau}$ is fine here.}\chisays{O.K. I changed all
%$\tilde{\tau}$ to $\tau$.}
\end{defn}

By the fusion rules for $L\left( \frac{1}{2},0\right) $, the $\tau
$-words for all irreducible $T$-submodules of $W$ are the same.
Thus, we can also define the {\it $\tau $-word of $W$} by
\[
\tau\left( W\right) =\tau(L(h_1,\cdots, h_n)),
\]
where $L(h_1,\cdots, h_n)$ is any irreducible $T$-submodule of
$W$.

The following proposition is an easy consequence of the fusion
rules (cf. \cite{DGH} and \cite{M2}).

\begin{prop}
\label{even}Let $C$ be an even code and let $W$ be an irreducible module of $%
M_{C}$. Then $\tau\left( W\right) $ is orthogonal to $C$.
\end{prop}

Now we shall give more details about the structure of the
irreducible module $W.$ The details can be found in \cite{M2}. Let
$\beta \in C^{\perp }:=\{\alpha\in \mathbb{Z}_2^n|\ \lan
\alpha,\gamma\ran=0 \ \text{ for all } \gamma\in C\}$ and
$C_{\beta }:=\{\alpha \in C|\ \mathrm{supp}\,\alpha \subseteq
\mathrm{supp}\,\beta \}$.

Let the group $\hat{C}=\left\{ \pm e^{k}|\ k\in C\right\} $ be a
central extension of $C$ by $\{\pm 1\}$ such that
$$e^he^k=(-1)^{\lan
h,k\ran}e^ke^h$$
 for any $h,k\in C$ and denote $\hat{C}_{\beta
}:=\left\{ \pm e^{k}|\ k\in C_{\beta}\right\}\subset \hat{C} $.
Let $H$ be a maximal self-orthogonal subcode of $C_{\beta }$. Then
$\hat{H}=\{\pm e^{\alpha }|\alpha \in H\}$ is a maximal abelian
subgroup of $\hat{C}_{\beta }$ (it is automatically normal since
it contains the commutator subgroup of $\hat{C}_{\beta }$). Take a linear character
$\chi : \hat{H}\rightarrow \{\pm 1\}$ with $\chi \left(
-e^{0}\right) =-1$ and define a 1-dimensional $\hat{H}$-module
$F_{\chi }$ by the action
\[
e^{\alpha }p=\chi \left( e^{\alpha }\right) p\qquad\text{ for }p\in F_{\chi },%
\text{ }\alpha \in H.
\]

We use ``$h_1\times h_2$'' to abbreviate a few of the well-known
fusion rules involving $L(\frac 1 2 , h_1)$ and $L(\frac 1 2,
h_2)$, i.e., $0\times h=h\times 0 =h$ for $h\in \{0,
\frac{1}2,\frac 1{16}\}$, $\frac{1}2\times \frac{1}2=0$ and
$\frac{1}2\times \frac{1}{16}=\frac{1}{16}\times
\frac{1}2=\frac{1}{16}$.

For any $h^{i}\in \{0,\frac{1}{2},\frac{1}{16}\},i=1,\cdots ,n$, with
${\tau}\left( \otimes_{i=1}^n L\left( \frac{1}{2},h^{i}\right)
\right) =\beta $, we define
\[
U=\left(\otimes _{i=1}^nL\left( \frac{1}{2},h^{i}\right)\right)
\otimes F_{\chi }.
\]
Then $U$ becomes an $M_{H}$-module with the vertex operator
defined by

\[
Y\left( \left( \otimes_{i=1}^n u^{i}\right) \otimes e^{\alpha
},z\right) =\left(\otimes
_{i=1}^{n}I^{\frac{a_{i}}{2},h^{i}}\left( u^{i},z\right)
\right)\otimes \chi \left( e^{\alpha }\right),
\]
where $u^{i}\in L(\frac{1}2,\frac{a_{i}}2)$, $\left( a_{1}, \dots,
a_{n}\right) \in H$, and $I^{\frac{a_{i}}{2},h^{i}}$ is 
a nonzero 
intertwining operator of type
\[
\left(
\begin{array}{ccc}
&L(\frac{1}2, \frac{a_i}2\times h^i)&\\ L(\frac{1}2,
\frac{a_i}2)&& L(\frac{1}2, h^i)
\end{array}
\right).
\]
We shall denote this $M_{H}$-module by $U\left( \left(
h^{i}\right) ,\chi \right) $ or $U\left( h^{i}\right) \otimes
F_{\chi }$.

%%% herehere

Let $\{\beta _{j}=\left( b_{j}^{i}\right) \}_{j=1}^{s}$ be a transversal of $%
H$ in $C$ and
\[
X=\bigoplus _{\beta _{j}\in C/H}\left\{U\left( h^{i}\times \frac{b_{j}^{i}}{2}%
\right) \otimes \left( e^{\beta _{j}}\otimes _{\hat{H}}F_{\chi
}\right) \right\},
\]
Note that $X$ does not
depend on the choice of the transversal of $H$ in $C$ and $X$ is an $M_{H}$%
-module.

The following results can  be found in Miyamoto \cite{M2}.

\begin{thm}
$X$ is an $M_{C}$-module with
\[
Y\left( u^{\gamma }\otimes e^{\gamma },z\right) =\left(
\otimes_{i=1}^n I\left( u^{i },z\right) \right) \otimes e^{\gamma
}
\]
for any $\gamma \in C$ and $u^{\gamma}=\otimes_{i=1}^n u^i\in
M_\gamma$.   We shall denote  $X$ by $\mathrm{Ind}_{H}^{C}U(
(h^i),\chi)$.
\end{thm}

\medskip

\begin{thm}\label{miya}
For any irreducible $M_{C}$-module $W$, there is a pair $\left(
\left( h^{i}\right) ,\chi \right) $ such that
\[
W\cong \mathrm{Ind}_{H}^{C}\left( U\left( \left( h^{i}\right)
,\chi \right) \right),
\]
where $\tau\left( W\right) =\tau\left( L\left( h^{1}, \cdots,
h^n\right) \right) =\beta $, $H$ is a maximal self-orthogonal
subcode of $C_{\beta}=\{\alpha \in C|\,\mathrm{supp}\,\alpha
\subseteq \mathrm{supp}\,\beta \}$ and $\chi $ is a linear character
of $\hat{H}$. Moreover, the structure of the $M_{C}$-module $W$ is
uniquely determined by an irreducible $M_{H}$-submodule of $W$.
\end{thm}

\medskip

Next we shall give a description of all irreducible $M_C$-modules
by using some binary words. Let $C$ be an even code of length $n.$
For a  given $\beta\in C^\perp$ and $\gamma\in \Z_2^n$, we define
\[
h_{\beta, \gamma}=(h_{\beta, \gamma}^1,\dots, h_{\beta,
\gamma}^n)\in \{0, \frac{1}2,\frac{1}{16}\}^n
\]
such that
\begin{equation*}
h_{\beta, \gamma}^i=
\begin{cases}
\frac{1}{16}& \text{ if }\quad  \beta_i=1,\\
\frac{\gamma_i}2 & \text{ if } \quad\beta_i=0.
\end{cases}
\end{equation*}
Denote $U(h_{_{\beta, \gamma}})=U(h_{\beta, \gamma}^1,\dots,
h_{\beta, \gamma}^n)=L(h_{\beta, \gamma}^1,\cdots,h_{\beta,
\gamma}^n)$. Fix a maximal self-orthogonal subcode $H_\beta$ of
the code $C_\beta=\{ \alpha\in C|\, \mathrm{supp}\,\alpha \subset
\mathrm{supp}\,\beta\}$ and define a character
$\chi_{\gamma}:\hat{H_\beta}\to \mathbb{C}$   of the abelian group
$\hat{H}_\beta$ by
\[
\chi_{\gamma}(-e^0)=-1 \quad \text{ and }\quad
\chi_{\gamma}(e^{\alpha})=(-1)^{\lan \alpha, \gamma\ran} \quad
\text{ for }\alpha\in H_\beta.
\]
Then $(\beta,\gamma)$ determines an irreducible $M_C$-module
\[
M_C(\beta,\gamma)=\mathrm{Ind}_{H_\beta}^C U(h_{\beta,
\gamma}^1,\dots, h_{\beta, \gamma}^n)\otimes F_{\chi_\gamma}.
\]
When there is no confusion, we shall simply denote
$M_C(\be,\gamma)$ by $M(\be, \gamma)$.

\begin{lem}
The definition of $M(\beta, \gamma)$ is independent of the choice
of the self-orthogonal subcode $H_\beta$ of $C_\beta$. 
\end{lem}

\pf Let $H$ be another maximal self-orthogonal subcode of
$C_\beta$ and let $\psi_\gamma: \hat{H} \to \mathbb{C}$ be a
character of $\hat{H}$ such that $\psi_\gamma(e^\xi)=(-1)^{\lan
\xi, \gamma\ran}$ and $\psi_\gamma(-e^0)=-1$. Then we can
construct another $M_C$-module
\[
\mathrm{Ind}_{H}^C U(h_{\beta, \gamma}^1,\dots, h_{\beta,
\gamma}^n)\otimes F_{\psi_\gamma}.
\]
By Miyamoto's Theorem\,(Theorem \ref{miya}), the structure of this
module is uniquely determined by the structure of the $M_H$
submodule $U(h_{\beta, \gamma}^1,\dots, h_{\beta,
\gamma}^n)\otimes F_{\psi_\gamma} $. Thus,
\[
\mathrm{Ind}_{H}^C U(h_{\beta, \gamma}^1,\dots, h_{\beta,
\gamma}^n)\otimes F_{\psi_\gamma}\cong  \mathrm{Ind}_{H_\beta}^C
U(h_{\beta, \gamma}^1,\dots, h_{\beta, \gamma}^n)\otimes
F_{\chi_\gamma}
\]
if and only if $ \mathrm{Ind}_{H_\beta}^C U(h_{\beta,
\gamma}^1,\dots, h_{\beta, \gamma}^n)\otimes F_{\chi_\gamma}$
contains an $M_H$-submodule isomorphic to $U(h_{\beta,
\gamma}^1,\dots, h_{\beta, \gamma}^n)\otimes F_{\psi_\gamma}. $ It
is equivalent to the fact that $
 \lan
\mathrm{Res}_{\hat{H}}\mathrm{Ind}_{\hat{H}_{\beta}}^{\hat{C_\beta}}{\chi_\gamma},
 \psi_\gamma \ran \neq 0,$ where
 $\lan
\mathrm{Res}_{\hat{H}}\mathrm{Ind}_{\hat{H}_{\beta}}^{\hat{C_\beta}}{\chi_\gamma},
 \psi_\gamma \ran$ denotes  the multiplicity of the character
 $\psi_\gamma$ in
 $\mathrm{Res}_{\hat{H}}\mathrm{Ind}_{\hat{H}_{\beta}}^{\hat{C_\beta}}{\chi_\gamma}$.

On the other hand,
\[
\mathrm{Res}_{\hat{H}}
\mathrm{Ind}_{\hat{H}_{\beta}}^{\widehat{H_{\beta}+H}}
{F_{\chi_\gamma}}\cong
 \bigoplus_{\alpha\in (H+H_\beta)/H_{\beta}}
e^\alpha\otimes F_{\chi_\gamma}\cong \bigoplus_{\alpha\in H/H\cap
H_{\beta}} e^\alpha\otimes F_{\chi_\gamma}.
\]

Let
\[
w_\psi= \frac{1}{\abs{H\cap H_\beta}}\sum_{\alpha \in H}
\psi_\gamma(\alpha) e^\alpha \otimes v,
\]
where $v\in F_{\chi_\gamma}$. Then $w_\psi \in
\mathrm{Ind}_{\hat{H}_{\beta}}^{\widehat{H_{\beta}+H}}
F_{\chi_\gamma}$ and for any $x\in H$,
\begin{align*}
e^x\cdot w_\psi &=\frac{1}{\abs{H\cap H_\beta}}\sum_{\alpha \in H}
(-1)^{\lan \gamma,\alpha\ran} e^x e^\alpha \otimes v \\
&=(-1)^{\lan \gamma,x\ran}\frac{1}{\abs{H\cap
H_\beta}}\sum_{\alpha \in H}
(-1)^{\lan \gamma,\alpha+x\ran} e^{x+ \alpha} \otimes v\\
&=(-1)^{\lan \gamma,x\ran} w_\psi=\psi_\gamma(e^x)w_\psi.
\end{align*}
Hence $\mathbb{C} w_\psi$ affords the $\hat{H}$-character
$\psi_\gamma$ inside
$\mathrm{Ind}_{\hat{H}_{\beta}}^{\widehat{H_{\beta}+H}} F_
{\chi_\gamma} \subset
\mathrm{Ind}_{\hat{H}_{\beta}}^{\hat{C_\beta}} F_{\chi_\gamma}$
and
\begin{align*}
 \lan
\mathrm{Res}_{\hat{H}}\mathrm{Ind}_{\hat{H}_{\beta}}^{\hat{C_\beta}}{\chi_\gamma},
 \psi_\gamma \ran \neq 0
\end{align*}
as desired. \qed

\begin{lem}\label{uniq}
Let $\beta_1, \beta_2\in C^\perp$ and $\gamma_1,\gamma_2\in
\Z_2^n$. Let $H_\beta$ be  a maximal self-orthogonal subcode of
$C_\beta$ and
 let $$(H_\beta)^{\perp_{\beta}}:=\{\alpha\in \Z_2^n\,|\, \mathrm{supp}\,
\alpha \subset \mathrm{supp}\,\beta \text{ and } \lan \alpha,
\xi\ran=0\ \text{ for all } \xi\in H_\beta\}.$$ Then the
irreducible $M_C$-modules $M(\beta_1,\gamma_1)$ and $
M(\beta_2,\gamma_2)$ are isomorphic if and only if
\[
\beta_1=\beta_2\qquad \text{ and }\qquad \gamma_1+\gamma_2\in
C+(H_\beta)^{\perp_{\beta}}.
\]
\end{lem}

\pf By the definition of $M(\beta,\gamma)$, it is easy to see that
$M(\beta_1,\gamma_1)\cong M(\beta_2,\gamma_2)$ if $
\beta_1=\beta_2$ and $\gamma_1+\gamma_2\in C$. Moreover, if
$\beta_1=\beta_2=\be$  and $ \gamma_1+\gamma_2\in
(H_\beta)^{\perp_{\beta}}$, then
$h_{\beta_1,\gamma_1}=h_{\beta_2,\gamma_2}$ and
$\chi_{\gamma_1}=\chi_{\gamma_2}$ for any choice of $H_\beta$.
Thus, $M(\beta_1,\gamma_1)\cong M(\beta_2,\gamma_2)$ if
\[
\beta_1=\beta_2\qquad \text{ and }\qquad \gamma_1+\gamma_2\in
C+(H_\beta)^{\perp_{\beta}},
\]

Now suppose that $M(\beta_1,\gamma_1)\cong M(\beta_2,\gamma_2)$.
Then they have the same $\tau$-word and $\beta_1=\beta_2$. Let
$\beta := \beta_1=\beta_2$. Let $H_\beta$ be a maximal
self-orthogonal subcode of $C_\beta$.   Since
$M(\beta_1,\gamma_1)\cong M(\beta_2,\gamma_2)$,
$M(\beta_1,\gamma_1)$  contains the $M_{H_\beta}$-module
$U(h_{\beta,\gamma_2})\otimes \chi_{\gamma_2}$. Thus, there exists
an element $\delta \in C$ such that
\[
h_{\beta, \gamma_1} \times \frac{\delta}2 = h_{\beta, \gamma_2}
\qquad \text{ and } \qquad  e^\delta \otimes F_{\chi_{\gamma_1}}
\cong F_{\chi_{\gamma_2}}.
\]
Since $h_{\beta, \gamma_1} \times \frac{\delta}2 = h_{\beta,
\gamma_2}$, $\delta+\gamma_1+\gamma_2 \in\Z_2^\beta$, where
$\Z_2^\beta=\{ \alpha\in {\Z_2}^n|\, \mathrm{supp}\,\alpha \subset
\mathrm{supp}\,\beta\}$. Moreover, $e^\delta \otimes
F_{\chi_{\gamma_1}} \cong F_{\chi_{\gamma_2}}$ implies that
\[
(-1)^{\lan \delta +\gamma_1,\alpha\ran}=(-1)^{\lan
\gamma_2,\alpha\ran}\qquad \text{ for all } \alpha\in H_\beta.
\]
Therefore, $\delta+\gamma_1+\gamma_2\in {H_\beta}^\perp$ and we
have $ \delta+\gamma_1+\gamma_2 \in {H_\beta}^\perp\cap
\Z_2^\beta=(H_\beta)^{\perp_{\beta}}$ and $\gamma_1+\gamma_2 \in
C+(H_\beta)^{\perp_{\beta}}. $ \qed

\begin{lem}
The code $C+H_\beta^{\bot_\beta}$ is independent of the choice of
the self-orthogonal subcode $H_\beta$.
\end{lem}

\proof

Let $H$ be another maximal self-orthogonal subcode of $C_\beta$.
Then we have $|H|=|H_\beta|$. First we will consider the
intersection $H\cap H_\be$ of $H$ and $H_\be$ and its orthogonal
complement in $\Z_2^\be$.

\medskip

\noindent \textbf{Claim:} $(H\cap H_\be)^{\perp_\be}=
H_\beta^{\bot_\be}+H$.

It is easy to see that $H_\beta^{\bot_\be}$ and $H$ are both
contained in $(H\cap H_\be)^{\perp_\be}$. Hence we have $
H_\beta^{\bot_\be}+ H \subset (H\cap H_\be)^{\perp_\be}$. Now note
that $H_\be^{\bot_{\be}}\cap H= H_\be^{\bot_{\be}}\cap (C_\be \cap
H)= (H_\be^{\bot_{\be}}\cap C_\be)\cap H=H_\be\cap H$ and $\dim H=
\dim H_\be$. By computing the dimensions, we have
\[
\begin{split}
\dim (H_\beta^{\bot_\be}+ H)&= \dim H_\beta^{\bot_\be}
+\dim H - \dim (H_\beta^{\bot_\be}\cap H )\\
&= (|\be| -\dim H_\beta) + \dim H - \dim
(H_\beta\cap H)\\
&=|\be| - \dim (H_\beta\cap H)\\
& =\dim (H\cap H_\be)^{\perp_\be}.
\end{split}
\]
Hence we have  $ (H\cap H_\be)^{\perp_\be}= H_\beta^{\bot_\be}+
H$.

\medskip

By the claim, we have
$$(H\cap H_\be)^{\perp_\be}= H+H_\be^{\perp_\be}
\subset C+H_\be^{\perp_\be}. $$ Therefore, $ C+(H\cap
H_\be)^{\perp_\be} \subset C+H_\be^{\perp_\be}$. On the other
hand, $C+H_\be^{\bot_\be}$ is clearly contained in $C+(H\cap
H_\be)^{\perp_\be}$ and thus $
 C+H_\be^{\perp_\be} = C+(H\cap
H_\be)^{\perp_\be}$. Similarly, we also have $ C+H^{\bot_\be}=
C+(H\cap H_\be)^{\perp_\be}$ and hence $C+H_\be^{\perp_\be}=
C+H^{\bot_\be}$ as desired. \qed

\medskip

Next  we shall compute the fusion rules among some irreducible
$M_C$-modules. We recall a theorem proved by
Miyamoto\,\cite{M2,M3}. Let $C$ be an even linear code.

\begin{thm}
For any $\al\in \Z_2^n$, the  $M_C$-module $\displaystyle M(0,
\al)= M_{\al+C}= \oplus_{\delta\in \al+C} M_\delta $ is a simple
current module. Moreover,
\[
M_{\al+C}\times M(\beta, \gamma) = M(\beta, \al+\gamma)
\]
for any irreducible $M_C$-module $M(\beta,\gamma)$.
\end{thm}

\medskip

Now by using the associativity and commutativity of the fusion
rules, we also have the following Lemma.

\begin{lem}\label{trans}
Let $\be_1,\be_2\in C^\perp$ and $\gamma\in \Z_2^n$. Then $$ \dim
I_{M_C} \binom{ M(\be_1+\be_2, \gamma)}{ M(\be_1,0)\qquad
M(\be_1,0)}= \dim I_{M_C} \binom{ M(\be_1+\be_2,
\al_1+\al_2+\gamma)}{ M(\be_1,\al_1)\qquad \qquad \qquad
M(\be_1,\al_2)}$$ for any $\al_1,\al_2 \in \Z_2^n$.
\end{lem}

\pf For any $\gamma\in \Z_2^n$, let
\[
m_\gamma= \dim I_{M_C} \binom{ M(\be_1+\be_2, \gamma)}{
M(\be_1,0)\qquad M(\be_1,0)}.
\]
Then we have
\[
M(\be_1, 0)\times M(\be_2, 0)= \sum_{\gamma \in \Z_2^n/K}
m_\gamma\, M(\be_1+\be_2, \gamma),
\]
where $K={C+(H_{\be_1+\be_2})^{\perp_{\be_1+\be_2}}}$.

Since the fusion product is associative and commutative, we have
\[
\begin{split}
&M(\be_1, \al_1)\times M(\be_2, \al_2)\\
 = &\left [M(
0, \al_1)\times M(\be_1, 0)\right] \times \left [M(0,
\al_2)\times M(\be_2, 0)\right]\\
=& \left [M(0,\al_1)\times M(0,\al_2)\right] \times \left
[M(\be_1,
0)\times M(\be_2, 0)\right]\\
=& M(0, \al_1+\al_2)\times \left [M(\be_1, 0)\times M(\be_2,
0)\right]\\
= & M(0, \al_1+\al_2)\times \big(\sum_{\gamma \in \Z_2^n/K}
m_\gamma\,  M(\be_1+\be_2, \gamma)\big)\\
=& \sum_{\gamma \in \Z_2^n/K} m_\gamma\,  M(\be_1+\be_2, \al_1+
\al_2+\gamma).
\end{split}
\]
Hence, we also have
\[
m_\gamma = \dim I_{M_C} \binom{ M(\be_1+\be_2,
\al_1+\al_2+\gamma)}{ M(\be_1,\al_1)\qquad \qquad \qquad
M(\be_1,\al_2)}
\]
as desired. \qed
\medskip

For the later purpose we also need some facts about the Hamming
code VOA $M_{H_8}$ from \cite{M2,M3} (see also \cite{La}).

Let $H_8$ be the Hamming code $[8,4,4]$ code, i.e., the code
generated by the rows of
\[
\begin{pmatrix}
1111& 1111\\
1111& 0000\\
1100& 1100\\
1010& 1010
\end{pmatrix}.
\]
Let $\{e^{1},\cdots ,e^{8}\}$ be the standard frame for $M_{H_8}$.
Let $q^0=\1$ be the vacuum element of $L(\frac 1 2,0)$ and let
$q^1$ be a highest weight vector of $L(\frac 1 2,\frac 1 2)$ such
that $q^1_0q^1=\1$. For any $\al=(\al_1, \dots, \al_8) \in H_8$,
let
\[
q^\al=q^{\al_1}\otimes \cdots \otimes q^{\al_8}\in M_\al,
\]
where $q^{\al_k}$ is a norm 1 highest weight vector for the $k$-th
tensor factor with respect to the action of our $T_8$. Then
$q^\al$ is a highest weight vector in $M_\al$. Moreover, we have
\[
{q^\al}_1 q^\be=
\begin{cases}
2 \sum_{1}^8 \al_i e^i &\text{ if } \al=\be,\\
q^{\al+\be} & \text{ if } |\,\al\cap \be|=2,\\
0& \text{ otherwise,}
\end{cases}
\]
%\bsays{I am unfamiliar with $<\al, \be>=2$.  I guess this means
%the codewords (as subsets) intersect in a 2-set, so why not just
%write $\al \cap \be$ is a 2-set or $| \al \cap \be |=2$.  }
%\chisays{I changed "$<\al, \be>$" to $|\al\cap\be|$ but I think
%people in coding theory use $<\al, \be>$ more often.}
for any
$\al,\be\in H_8$ with $|\al| =|\be|=4$.

The following results are obtained in \cite{M2}.

\begin{lem}
\label{h8} Let $\nu _{i}$ be the binary word whose $i$-th entry is
$1$ and all other entries are $0$.  Define $\alpha _{i}:=\nu
_{1}+\nu _{i}$.

In the Hamming code VOA $ M_{H_{8}}$, there exist exactly three
Virasoro frames, namely,
\begin{equation*}
\{e^{1},\cdots ,e^{8}\},\left\{ d^{1},\cdots ,d^{8}\right\} ,\text{ and }%
\left\{ f^{1,}\cdots ,f^{8}\right\}
\end{equation*}
where
$$d^{i}=S^{\alpha ^{i}}=\frac{1}{8}(e^{1}+\cdots +e^{8})+\frac{1}{8}\sum_{\beta \in H_{8},\left| \beta \right| =4}(-1)^{\<\alpha
_{i},\beta \>}q^{\beta }\otimes e^{\beta },$$
%$$\alpha _{i}=\nu _{1}+\nu _{i},$$
$$f^{i}=S^{\nu _{i}}=\frac{1}{8}(e^{1}+\cdots
+e^{8})+\frac{1}{8}\sum_{\beta \in H_{8},\left| \beta \right|
=4}(-1)^{\<\nu _{i},\beta \>}q^{\beta }\otimes e^{\beta }$$.
%and
%$\nu _{i}$ is the binary word whose $i$-th entry is $1$ and all
%other entries are $0$.
\end{lem}

\begin{thm} Let $L$ be an irreducible
$M_{H_{8}}$-module with half-integral or integral weight. Then,
$L$ is isomorphic to one of the following:

\begin{enumerate}
\item  $M_{\nu _{1}+\nu _{i}+H_{8}}$ with respect to $\{e^{1},\cdots ,e^{8}\}
$ for all $i=1,\cdots ,8.$

\item  $M_{\nu _{i}+H_{8}}$ with respect to $\{e^{1},\cdots ,e^{8}\}$ for
all $i=1,\cdots ,8.$

\item  $M_{\nu _{i}+H_{8}}$ with respect to $\left\{ d^{1},\cdots
,d^{8}\right\} $ for all $i=1,\cdots ,8.$

\item  $M_{\nu _{i}+H_{8}}$ with respect to $\left\{ f^{1,}\cdots
,f^{8}\right\} $ for all $i=1,\cdots ,8.$
\end{enumerate}
Moreover, all modules in (3) and (4) are isomorphic to $\otimes
_{i=1}^{8}L(\frac{1}{2},\frac{1}{16})$ as $T_8$-modules.
\end{thm}

As a corollary, we have the following theorem. The proof can be
found in \cite{La} (see also \cite{M2,M3}).

%\begin{thm}\label{fh8}
%All irreducible $M_{H_{8}}$-modules  with half-integral or
%integral weight are simple current modules. Moreover,
%\[
%M(\be_1, \al_1)\times_{_{M_{H_8}}} M(\be_2, \al_2)= M(\be_1+
%\be_2, \al_1+\al_2)
%\]
%for any $\be_1= (0^8) $ or $(1^8)$ and $\be_2\in H_8$.
%\end{thm}

\begin{thm}\label{fh8}
For any $\be_1= (0^8) $ or $(1^8)$ and $\be_2\in H_8$, we have
\[
M(\be_1, \al_1)\times_{_{M_{H_8}}} M(\be_2, \al_2)= M(\be_1+
\be_2, \al_1+\al_2)
\]
Consequently, all irreducible $M_{H_{8}}$-modules with
half-integral or integral weight are simple current modules.
\end{thm}

\begin{rem}
For any $\al\in \Z_2^8/ H_8$, $\al$ uniquely determines a
character $\chi_{\al}\in \mathrm{Irr}\, H_8$ such that
$\chi_{\al}(\gamma) =(-1)^{\lan \al, \gamma\ran}$ for any $\gamma
\in H_8$. By using this identification, our module $M(\be, \al)$
actually corresponds to the class $[\be, \chi_\al]$ defined in
Section 5 of \cite{La}.
\end{rem}

\section{The moonshine vertex operator algebra  $V^{\natural}$}
\setcounter{equation}{0}

Let $V^{\natural}$ be the moonshine vertex operator algebra
\cite{FLM1}-\cite{FLM2}. The following theorem can be found in
\cite{DGH}.

\begin{thm}\label{t3.1}
There exists a VF in  $V^{\natural}$, called 
$F:=\{\omega_1,...,\omega_{48}\}$, so that the code $\sC:=C(F)$
associated to this VF has length $48$ and dimension $41$. The code
$\sD:=D(F)=\sC^{\perp}$ has generator matrix
$$\left(\begin{array}{ccc}
1111 1111 1111 1111 & 0000 0000 0000 0000 & 0000 0000 0000 0000 \\
0000 0000 0000 0000 & 1111 1111 1111 1111 & 0000 0000 0000 0000 \\
0000 0000 0000 0000 & 0000 0000 0000 0000 & 1111 1111 1111 1111 \\
0000 0000 1111 1111 & 0000 0000 1111 1111 & 0000 0000 1111 1111 \\
0000 1111 0000 1111 & 0000 1111 0000 1111 & 0000 1111 0000 1111 \\
0011 0011 0011 0011 & 0011 0011 0011 0011 & 0011 0011 0011 0011 \\
0101 0101 0101 0101 & 0101 0101 0101 0101 & 0101 0101 0101 0101
\end{array}\right)_{\textstyle .}$$
\end{thm}

\begin{rem} The weight enumerator of $\sD$  is
given by $$X^{48}+3X^{32}+120X^{24}+3X^{16}+1$$ and the minimal
weight of $\sC$ is $4$. Moreover, $\sD$ is self-orthogonal and
hence $\sD \subset \sC.$  (The codes $\sD$ and $\sC$ are denoted
by $S^\natural$ and $D^\natural$, respectively in \cite{M3}.)
\end{rem}

\begin{lem}\label{codec} The code $\sC$ in Theorem \ref{t3.1} is generated
by the weight 4 codewords.
\end{lem}

\pf First we note that the code $\sD=D(F)$ is generated by the
elements of the form
\[
(1^{16},0^{16},0^{16}),\ (0^{16},1^{16},0^{16}), \
(0^{16},0^{16},1^{16}) \quad \text{ and }\quad (\al, \al, \al) ,
\quad \al \in \mathrm{RM}(1,4),
\]
where $\mathrm{RM}(r,m)$ denote the $r$-th order Reed-Muller code
of length $2^m$ (cf. \cite{CS}).

Since $\mathrm{RM}(1,4)^\perp= \mathrm{RM}(2,4)$, we have
\[
\sC=\sD^\perp=\{ (\al, \be, \gamma)|\ \al+\be+\gamma\in
\mathrm{RM}(2,4),\  \al, \be, \gamma \text{ even }\}.
\]
Hence the code $\sC$ can be generated by the elements
\[
(\al, 0, 0),\ (0, \be, 0), (0, 0, \gamma), \quad \al, \be, \gamma
\text{ are generators of }\mathrm{RM}(2,4)
\]
and
\[
(\al, \be, 0),\ (\al, 0, \be), (0, \al, \be), \quad \al, \be
\text{ are even and } \al+\be \text{ is a generator of }
\mathrm{RM}(2,4).
\]
Note that the Reed Muller code $\mathrm{RM}(2,4)$ is of dimension
$11$ and is generated by the elements of the form
\[
(\al, 0),\quad  (0,\al), \quad \al\in H_8,
\]
and
\[
(1100\, 0000\, 1100\, 0000), \ (1010\, 0000\, 1010\, 0000), \
(1000\, 1000\, 1000\, 10000).
\]

\medskip

Since the Hamming code $H_8$ is generated by its weight $4$
elements, the codes $\mathrm{RM}(2,4)$ and  $\sC$ are  generated
by the  weight $4$ codewords also.  \qed

In the next theorem, see \ref{DGH}(d) for the meaning of
$(V^{\natural})^0$.

\begin{lem} The vertex operator subalgebra $(V^{\natural})^0$ is
isomorphic to $M_{\sC}$ and  is uniquely determined by the set of
weight 4 codewords of $\sC.$
\end{lem}

\pf By the uniqueness of the code VOA, $(V^{\natural})^0$ and
$M_{\sC}$ are isomorphic. Since $\sC$ is generated by the weight 4
codewords of $\sC,$ the vertex operator algebra structure of
$(V^{\natural})^0$ is uniquely determined by the generators of the
group $\sC.$ \qed

We now determine the irreducible modules and the fusion rules for
the code VOA $M_{\sC}$.

\begin{rem}
In the next result,
 $\mathrm{RM}(r,m)$ denote the $r$-th order Reed-Muller code
of length $2^m$ (cf. \cite{CS}). Note that the Reed Muller codes
are nested in the sense that  $\mathrm{RM}(r+1, m)\supset
\mathrm{RM}(r,m)$ and $\mathrm{RM}(r+1, m+1)\supset
\mathrm{RM}(r,m)\oplus \mathrm{RM}(r,m)$, where the direct sum
corresponds to a partition of indices by an affine hyperplane and
its complement.
\end{rem}

The following properties of the code $\sC$ can be derived easily
from the definition.

\begin{prop}
Let $\sD$ and $\sC$ be defined as above. For any $\be\in \sD$,
denote $$\sC_\beta :=\{ \alpha\in \sC|\ \mathrm{supp}\,\alpha
\subset \mathrm{supp}\,\beta\}.$$
\begin{enumerate}

\item  If $|\beta|=16$, then $\sC_\beta\cong \mathrm{RM}(2,4)$.

\item If $|\beta|=24$, then
$ \sC_\beta\cong \{ (\al, \gamma, \delta)|\ \al+\gamma+\delta \in
H_8 \text{ and } \al,\gamma, \delta \text{ even}\}.  $

\item If  $|\beta|=32$, then $ \sC_\beta\cong \mathrm{RM}(3,5).$

\item If $|\beta|=48$, then $ \sC_\beta = \sC$.
\end{enumerate}
Note that the Hamming code $H_8\cong \mathrm{RM}(1,3)$. Hence, for
$\be \ne 0$, $\sC_\be$  contains a self-dual subcode which is
isomorphic to a direct sum of $|\be|/8$ copies of the Hamming code
$H_8$.
\end{prop}

\pf Let $\be\in \sD$ and  $n=|\be|$,  the weight of $\beta$.  Let
$p_\be: \Z_2^{48}\to \Z_2^n$  be the natural projection of
$\Z_2^{48}$ to the support of $\beta$. Since $\sC= \sD^\perp$, it
is easy to see that $\sC_\be\cong {p_\be(\sD)}^\perp$.

\medskip

\noindent \textbf{Case 1.} $|\beta|=16$. In this case,  $
p_\be(\sD)$ is generated by the codewords $(1^{16})$, $(0^8\,1^8),
(0^4\ 1^4)^2$, $(0^2\ 1^2)^4$ and $(0\ 1)^8$ and is isomorphic to
the Reed Muller code $\mathrm{RM}(1,4)$. Since
$$\mathrm{RM}(r,m)^\perp \cong \mathrm{RM}(m-r-1,m)$$
 for any
$1\leq r\leq m$ we have $\sC_\be\cong \mathrm{RM}(2,4)$ as
desired.

\bigskip

\noindent \textbf{Case 2.} $|\beta|=24$. In this case,
$p_\be(\sD)$ is of dimension $6$ and is isomorphic to a code
generated by
\[
(1^8\, 0^8\, 0^8),\  (0^8\, 1^8\, 0^8),\ (0^8, 0^8, 1^8)\  \text{
and }\ (\al, \al,\al),\ \al\in H_8.
\]
Hence  $\sC_\be\cong \{ (\al, \gamma, \delta)|\ \al+\gamma+\delta
\in H_8 \text{ and } \al,\gamma, \delta \text{ even}\}$.

\bigskip

\noindent \textbf{Case 3.} $|\beta|=32$.  $p_\be(\sD)\cong
\mathrm{RM}(1,5)$ and hence $\sC_\be\cong \mathrm{RM}(3,5)$.

\bigskip

\noindent \textbf{Case 4.} $|\beta|=48$. It is clear that
$\sC_\be=\sC$ in this case. \qed

\medskip

Now by using Lemma \ref{uniq}, we have the following theorem.
\begin{thm}\label{indeed}
Let $\sD$ and $\sC$ be defined as above. Then
\[
\{ M_\sC(\be, \gamma)\, |\ \be\in \sD \text{ and } \gamma \in
\Z_2^{48}/ \sC\}
\]
is the set of all inequivalent irreducible modules for $M_\sC$.
\end{thm}

\pf By the previous proposition, we can choose $H_\be$ such that
it is a direct sum of $|\be|/8$ copies of the Hamming code $H_8$.
In this case, $(H_\beta)^{\perp_{\beta}}= H_\be$ and we have $\sC=
\sC+(H_\beta)^{\perp_{\beta}}$. Hence $ \{ M_\sC(\be, \gamma)\, |\
\be\in \sD \text{ and } \gamma \in \Z_2^{48}/ \sC\} $ is the set
of all inequivalent irreducible modules for $M_\sC$ by Lemma
\ref{uniq} . \qed

\medskip

Next  we shall compute the fusion rules among irreducible
$M_\sC$-modules. The main tool is the representation theory of the
Hamming code VOA $M_{H_8}$ given in Section 4. First we recall the
following theorem from \cite{DL}.  
%\bsays{Is this true for
%generating sets?  Say if X generates W1 and Y generated W2 and if
%$I(x,z)y=0$ for all x in X and y in Y, then is $I(\cdot , z)=0$? }
%\chisays{ Yes, it is true. It can be proved easily using
%associativity.}\bsays{It would seem worthwhile to prove this natural
%extension.  How would you two feel about doing so here?  }
%\chosays{I dont think we need this. It is a trivial result.}
%\chisays{I also feel it is a trivial result}

\begin{thm}\label{dandl} Let $W^{1},W^{2}$ and $W^{3}$ be $V$-modules and let
$I$ be an intertwining operator of type
\[
\left(
\begin{array}{ccc}
& W^{3} &  \\ W^{1} &  & W^{2}
\end{array}
\right) .
\]
Assume that $W^{1}$ and $W^{2}$ have no proper submodules
containing $v^{1}$ and $v^{2}$, respectively. Then $I\left(
v^{1},z\right) v^{2}=0\text{ implies }I\left( \cdot ,z\right) =0.$
\end{thm}

\begin{lem}\label{ub} For any $\be_1, \be_2,
\be_3\in \sD$ and $\al_1, \al_2, \al_3\in \Z_2^{48}$, we have
\[
\dim I_{M_\sC} \binom{ M(\be_3, \al_3)}{M(\be_1, \al_1)\qquad
M(\be_2, \al_2)} \leq 1
\]
and
\[
\dim I_{M_\sC} \binom{ M(\be_3, \al_3)}{M(\be_1, \al_1)\qquad
M(\be_2, \al_2)}=0
\]
unless $\be_3=\be_1+\be_2$ and $ \al_3 \equiv \al_1+\al_2  \ \mod
\ \sC$.
\end{lem}

\pf Recall that $\dim \sD=7$ and the weight enumerator of $\sD$ is
$X^{48}+3X^{32}+120X^{24}+3X^{16}+1.$

Without loss, we may assume that $\be_3=\be_1+\be_2$; otherwise,
\[
\dim I_{M_\sC} \binom{ M(\be_3, \al_3)}{M(\be_1, \al_1)\qquad
M(\be_2, \al_2)}=0.
\]

Let $\bar{\be_1}= (1^{48})+ \be_1$. Then $\bar{\be_1}$ is also in
$\sD$. Thus, there exist self-orthogonal codes $H_{\be_1}$ and
$H_{\bar{\be_1}}$ of $\sC$ such that both $H_{\be_1}$ and
$H_{\bar{\be_1}}$ are direct sums of Hamming $[8,4,4]$ codes. Let
$E= H_{\be_1}\oplus H_{\bar{\be_1}}\cong {H_8}^{\oplus 6}$. If the
weight of $\be_2$ is a multiple of $16$ (i.e., $0,16, 32,$ or
$48$), then $|\supp\, \be_1 \cap \supp\, \be_2|$ is a multiple of
$8$. In this case, it is possible to find maximal self-orthogonal
subcodes $H_{\be_2}$ of $\sC_{\be_2}$ and $H_{\be_1+\be_2}$ of
$\sC_{\be_1+\be_2}$ such that  $H_{\be_2}$ and  $H_{\be_1+\be_2}$
are isomorphic to direct sums of Hamming codes and are both
contained in $E$. Then as an $M_E$-module,
\[
M_{\sC}(\be_i, \al_i)= \bigoplus_{\delta\in \sC/ E} M_{E}(\be_i,
\al_i+\delta).
\]
Note that $H_{\be_i}\subset E$  for any $i=1,2,3$ and hence
$M_{\sC}(\be_i, \al_i)$ is a direct sum of inequivalent
irreducible $M_E$-modules. Thus by Theorem \ref{fh8} and
\ref{dandl}, we have
\[
\dim I_{M_\sC} \binom{ M(\be_1+\be_2, \al_3)}{M(\be_1,
\al_1)\qquad M(\be_2, \al_2)} \leq \dim I_{M_E} \binom{
M_E(\be_1+\be_2, \al_3)}{M_E(\be_1, \al_1)\qquad M_E(\be_2,
\al_2)}\leq 1
\]
and
\[
\dim I_{M_\sC} \binom{ M(\be_1+\be_2, \al_3)}{M(\be_1,
\al_1)\qquad M(\be_2, \al_2)}=0
\]
unless $\al_3=\al_1+\al_2$.

\medskip

Finally, we shall treat the case for which all $\be_1$, $\be_2$
and $\be_1+\be_2$ are of weight $24$. For simplicity, we may
assume that $\be_1=(1^8\, 0^8\, 1^8\, 0^8\,1^8\, 0^8)$ and $\be_2=
(1^4\, 0^4\, \dots\, 1^4\, 0^4)$. Then $\be_3=\be_1+\be_2= (0^4\,
1^8\, 0^8\, 1^8\, 0^8\,1^8\,0^4)$.  In this case, we have $E=
H_{\be_1}\oplus H_{\bar{\be}_1}\cong {H_8}^{\oplus 6}$,
$H_{\be_2}\cong H_8\oplus H_8 \oplus H_8$ and
$H_{\be_1+\be_2}\cong H_8\oplus H_8 \oplus H_8$.

\medskip

Note that $E+H_{\be_1+ \be_2}=E+H_{\be_2}$ in this case. Moreover,
we have
\[
E_{\be_2}=\{ \alpha\in E|\, \mathrm{supp}\,\alpha \subset
\mathrm{supp}\,\beta_2\}= E\cap H_{\be_2}\quad \text{ and }
E_{\be_3}=E_{\be_1+\be_2}= E\cap H_{\be_1+ \be_2}.
\]
Let $\mathcal{H} := E+ H_{\be_2}= E+H_{\be_1+\be_2}$. Then the
$M_{\sC}$-module $M_\sC(\be_i, \al_i), i=2,3,$ can be decomposed
as
\[
M_{\sC}(\be_i, \al_i)= \bigoplus_{\delta\in \sC/ \mathcal{H}}
M_{\mathcal{H}}(\be_i, \al_i+\delta).
\]

\medskip

\noindent \textbf{Claim: } $M_{\mathcal{H}}(\be_i, \al_i+\delta)$
is irreducible as an $M_E$-module for any $\delta\in \sC/\sH$.

\medskip

Proof. Let $W=M_{\mathcal{H}}(\be_i, \al_i+\delta)$ and
$\sH_{\be_i}= \{ \al\in \sH|\, \supp \, \al \subset \supp \,
\be_i\}$. Then $H_{\be_i}$ is a maximal self-orthogonal subcode of
$\sH_{\be_i}$. Let $U( \mathbf{h})\otimes F_{\chi}$ be an
irreducible $M_{H_{\be_i}}$-submodule of $W$. Then
$$W=\mathrm{Ind}_{H_{\be_i}}^\sH U( \mathbf{h})\otimes F_{\chi}=
\bigoplus_{\delta\in \sH/ H_{\be_i}} U( \mathbf{h}\times
\frac{\delta}2)\otimes( e^\delta\otimes F_{\chi}). $$

Since $E_{\be_i}= E\cap H_{\be_i}\subset H_{\be_i}$,
$U(\mathbf{h})\otimes F_{\chi}$ is also an irreducible
$M_{E_{\be_i}}$-module. Hence,
\[
W'= \mathrm{Ind}_{E_{\be_i}}^E U( \mathbf{h})\otimes F_{\chi}=
\bigoplus_{\delta\in E/ E_{\be_i}} U( \mathbf{h}\times
\frac{\delta}2)\otimes( e^\delta\otimes F_{\chi})
\]
is an irreducible $M_E$-submodule of $W$.  Note that
\[
|\sH/ H_{\be_i}|= | (E+H_{\be_i})/ H_{\be_i}|=| E/ (E\cap
H_{\be_i})|= |E/ E_{\be_i}|.
\]
Therefore, we have $W=W'$ and $W$ is an irreducible $M_E$-module.

Now, by Theorem \ref{fh8} and \ref{dandl}, we have
\[
\begin{split}
& \dim I_{M_\sC} \binom{ M(\be_1+\be_2, \al_3)}{M(\be_1,
\al_1)\qquad M(\be_2, \al_2)} \\
\leq & \dim I_{M_\sH} \binom{ M_\sH(\be_1+\be_2,
\al_3)}{M_\sH(\be_1, \al_1)\qquad M_\sH(\be_2, \al_2)}\\
\leq &\dim I_{M_E} \binom{ M_E(\be_1+\be_2, \al_3)}{M_E(\be_1,
\al_1)\qquad M_E(\be_2, \al_2)}\leq 1
\end{split}
\]
and
\[
\dim I_{M_\sC} \binom{ M(\be_1+\be_2, \al_3)}{M(\be_1,
\al_1)\qquad M(\be_2, \al_2)}=0
\]
unless $\al_3=\al_1+\al_2$. Note that $M_{\sH}(\be_2,\al_2)=
M_E(\be_2,\al_2)$ and $M_{\sH}(\be_1+\be_2,\al_1+\al_2)=
M_E(\be_1+\be_2,\al_1+\al_2)$ as $M_E$-modules. \qed

\medskip

\begin{thm}\label{fusion}
The fusion rules among irreducible $M_\sC$ modules are given by
\[
M(\be_1, \al_1) \times M(\be_2, \al_2) = M(\be_1+\be_2,
\al_1+\al_2),
\]
where $\be_1, \be_2\in \sD$ and $\al_1, \al_2\in \Z_2^{48}/\sC$.
In particular, each irreducible $M_{\sC}$-module is a simple
current.
\end{thm}

\pf By Lemma \ref{trans} and \ref{ub},  it remains to show that
\[
I_{M_\sC} \binom{ M(\be_1+\be_2, \al_1+\al_2)}{M(\be_1,
\al_1)\qquad M(\be_2, \al_2)}\neq 0,
\]
for some $\al_1, \al_2\in \Z_2^{48}$.  Nevertheless, such kind of
intertwining operators does exist and can be realized inside the
Leech lattice VOA $V_\Lambda$. In fact, there exists a Virasoro
frame of $V_\Lambda$ such that $V_\Lambda$ can be decomposed as
\[
V_{\Lambda}\cong \bigoplus_{\be\in \sD} M_\sC(\be, \gamma_{\be}),
\quad \text{ for some } \gamma_{\be}\in \Z_2^{48}/ \sC.
\]
We shall refer to \cite{DGH} or \cite{M3} for details. \qed

\section{Proof of the main theorems}
\setcounter{equation}{0}

We first prove Theorem \ref{mt1}. So we assume that (1)  $V$ is a
vertex operator algebra satisfying conditions (a)-(c), (2) $V_2$
is isomorphic to the Griess algebra, (3) $V$ is $C_2$-cofinite.

\begin{lem}\label{l4.1} $V$ is truncated below 0 and $V_0=\C \1.$
\end{lem}

\pf   First we prove that $V_n=0$ if $n$ is negative. If this is
not true,  take the smallest $n$ such that $V_n\ne 0.$ Then each
$0\ne v\in V_n$ generates a highest weight module for the Lie
algebra $\C L(1)\oplus \C L(0)\oplus \C L(-1)$ (which is
isomorphic to $sl(2,\C).$) According to the structure of the
highest weight modules for $sl(2,\C)$ we know that $L(-1)^iv\ne 0$
for $i=0,...,-2n.$   Since $n$ is less than or equal to $-1$ we
see that $L(-1)^{-n+1}v\ne 0.$ Since the weight of $L(-1)^{-n+1}v$
is 1, we immediately have a contradiction as $V_1=0$ by
assumption.

We now prove that $V_0$ is one dimensional. Note that
$L(-1)V_0=0.$ So each nonzero vector $v\in V_0$ is a vacuum-like
vector \cite{Li}.
As a result, we have a $V$-module
isomorphism  $f_v: V \to V$ by sending $u$ to $u_{-1}v$ for $u\in
V$ \cite{Li}.   By Schur's lemma, $f_v$ must be a multiple of the
identity map. As a result, $f_v(\1)=v$ is a multiple of the
vacuum. This shows that $V_0$ is spanned by the vacuum. \qed

\bigskip

\begin{lem}\label{l4.2r} $V$ is a holomorphic vertex operator algebra.
\end{lem}

\pf It is proved in \cite{DLM2} that if $U$ is a vertex operator
algebra such that $U=\oplus_{n\geq 0}U_n$ with $U_0$ being
1-dimensional and $U_1=0$ and that $U$ is the only irreducible
ordinary module for itself then any ordinary module is completely
reducible. So any ordinary $V$-module is a direct sum of copies of
$V.$ Since $V$ is $C_2$-cofinite, any submodule generated by a
single vector in any admissible module is an ordinary module (see
\cite{ABD}). This shows that any admissible $V$-module is
completely reducible. That is, $V$ is rational. This together with
condition (a) gives conclusion that $V$ is holomorphic. \qed

\bigskip

\begin{lem}\label{lj} The $q$-dimension  $ch_q V=q^{-1}\sum_{n\geq 0}(\dim
V_n)q^n$ of $V$  is $J(q).$
\end{lem}

\pf Since $V$ is holomorphic and $C_2$-cofinite, by the modular
invariance result in \cite{Z}, $ch_q V$ is a modular function on
the full modular group, and thus equal to $J(q)$ by noting that
$V_0=\C {\bf 1}$ and $V_1=0.$ \qed
\bigskip

Since $V$ is irreducible and $V_0/L(1)V_1$ is one dimensional,
there is a unique nondegenerate symmetric invariant bilinear form
$(\cdot,\cdot)$ on $V$ such that $(\1,\1)=1$ (see \cite{Li}). That
is,
$$(Y(u,z)v,w)=(-z^{-2})^{\wt u} (v,Y(e^{zL(1)}u,z^{-1})w)$$
for homogeneous $u\in V.$ In particular, the restriction of
$(\cdot,\cdot)$ to each $V_n$ is nondegenerate. As a result,
$(\cdot,\cdot)$ defines a nondegenerate symmetric invariant
bilinear form on the Griess algebra $V_2$ such that $(u,v)=u_3v$
for $u,v\in V_2.$

From now on we will fix the vectors $\{\omega_1,...,\omega_{48}\}$
of $V_2$ given in Theorem \ref{t3.1}. Since we only assume that
$V_2$ is isomorphic to the Griess algebra we do not know if the
bilinear form $(u,v)=u_3v$ defined on $V_2$ is the same as the
bilinear form defined on $V^{\natural}_2$ using the same formula.
So it is not clear that $\{\omega_1,...,\omega_{48}\}$ forms a VF
in $V.$

Since $V_2$ is a simple commutative nonassociative algebra, we
need a result on the bilinear forms over
 a finite dimensional simple commutative nonassociative
algebra $B.$ A bilinear form $(\cdot,\cdot)$ on $B$ is called {\it
invariant} if $(ab,c)=(b,ac),$ for all $a, b, c \in B$.   The next
result applies to any finite dimensional simple algebra.

\begin{lem}\label{lb} The space of nondegenerate symmetric invariant
bilinear forms on $B$ is at most one-dimensional.
\end{lem}

\pf Let $(\cdot,\cdot)$ and $\<\cdot,\cdot\>$ be two nondegenerate
symmetric invariant bilinear forms on $B.$ Then there is a linear
isomorphism $f: B\to B$ such that $(u,v)=\<f(u),v\>$ for all
$u,v\in V.$ For any $a\in B$ we have
$$\<f(au),v\>=(au,v)=(u,av)=\<f(u),av\>=\<af(u),v\>.$$
That is, $f(au)=af(u).$ Let $B_{\lambda}$ be the eigenspace of $f$
with eigenvalue $\l\ne 0$ Then $B_{\lambda}$ is an ideal of $B.$
This shows that $B=B_{\lambda}.$ So $f=\lambda \ id_B.$ As a
 result, $(\cdot,\cdot)=\l\<\cdot,\cdot\>,$ as desired.
\qed

\begin{lem}\label{ladd} Each $\omega_i$ is a Virasoro vector with central
charge $\frac 1 2$ and  for all $m, n$,
$$[L^i(m),L^j(n)]=0$$
if
$i\ne j$ where $Y(\omega_i,z)=\sum_{n\in\Z}L^i(n)z^{-n-2}.$
\end{lem}

\pf  We first prove that each $\omega_i$ is a Virasoro vector of
central charge $\frac{1}{2}.$ That is, the component operators
$L^i(n)$ of $Y(\omega_i,z)=\sum_{n\in\Z}L^i(n)z^{-n-2}$ satisfies
the Virasoro algebra relation with central charge $\frac{1}{2}.$

Clearly $\omega_i\cdot\omega_i=L^i(0)\omega_i=2\omega_i$ by the
product in $B.$ So, $\omega_i$ is a Virasoro vector with central
charge $c_i$ defined by  $c_i{\bf 1} =2L^i(2)\omega_i.$ Note that
$L^i(0)$ is semisimple on $V_2$ and the eigenvalues of $L^i(0)$
are $2,0,\frac{1}{2}$ and $\frac{1}{16}$ (see \cite{DGH}). Since
the bilinear form is invariant, we see that the eigenspaces with
different eigenvalues are orthogonal. So the restriction of the
bilinear form to each eigenspace is nondegenerate. It is known
from \cite{DGH} that the eigenspace with eigenvalue $2$ is one
dimensional and is spanned by $\omega_i.$ As a result,
$L^i(0)\omega_i$ is nonzero and $c_i\ne 0.$ We must prove that
$c_i=\frac{1}{2}.$

Recall from \cite{DM3}  that the Griess algebra is a simple
commutative nonassociative algebra. Let $\<\cdot,\cdot\>$ be the
bilinear from defined on $V_2^{\natural}$ and $(\cdot,\cdot)$ be
the bilinear form defined on $V_2.$ By Lemma \ref{lb},
$(\cdot,\cdot)$ is a multiple of $\<\cdot,\cdot\>.$ Note that
$\<\omega,\omega\> =(\omega,\omega)=12.$ We conclude that these
two bilinear forms are exactly the same. So
$(\omega_i,\omega_i)=\<\omega_i,\omega_i\>=\frac{1}{4}.$ That is
$c_i=\frac{1}{2}.$

Let $i\ne j.$ Since $(\omega_i,\omega_j)=0$ and
$L^i(0)\omega_j=0$, we see immediately that  $[L^i(m),L^j(n)]=0$
for all $m,n\in \Z.$ \qed

\begin{thm}\label{t4.2} The  $\{\omega_1,...,\omega_{48}\}$ forms a VF in $V$
and $V$ is a FVOA.
\end{thm}

\pf We only need to prove that vertex operator subalgebra
$\<\omega_i\>$ generated by $\omega_i$ is isomorphic to
$L(\frac{1}{2},0)$ for the Virasoro algebra $Vir_i$ generated by
$L^i(m)$ for $m\in \Z.$ It is clear that $\<\omega_i\>$ is a
highest weight module with highest weight $0$ for $Vir_i.$ Then
there are two possibilities. Either $\<\omega_i\>$ is the Verma
module modulo the submodule generated by $L^i(-1)\1$ or
 $\<\omega_i\>$  is isomorphic to $L(\frac{1}{2},0)$, according to the structure theory of highest weight modules for the Virasoro algebra with central charge $\frac 1 2$ \cite{FF}.    We now assume that
the first possibility happens.

In this case the $q$-character of $\<\omega_i\>$ is equal to
$$ch_q\<\omega_i\>=q^{-1/48}\frac{1}{\prod_{n\geq 2}(1-q^n)}.$$
Let $U$ be the vertex operator subalgebra of $V$ generated by
$\omega_j$ for $j=1,...,48.$ Then we have
$$U=\<\omega_1\>\otimes \cdots \otimes \<\omega_{48}\>$$
is a tensor product. Let $ f(q^{\frac 1 n }),g(q^{\frac 1 n })\in
\R[[q^{1/n},q^{-1/n}]]$ for some positive integer $n.$ We write $
f(q^{\frac 1 n })\leq g(q^{\frac 1 n })$ if the coefficient of
$q^m$ in $ f(q^{\frac 1 n })$ is less than or equal to that in
$g(q^{\frac 1 n })$ for all $m.$ It is well known  that the
$q$-character of $L(\frac{1}{2},0)$ is equal to
$$\frac{1}{2}q^{-1/48}\left(\prod_{n\geq 0}(1+q^{n+\frac 1 2})+\prod_{n\geq 0}(1-q^{n+\frac 1 2})\right)$$
(cf. \cite{KR}). Thus we have
\begin{equation}
\label{ine} ch_qU\geq  f(q^{\frac 1 2 })
\end{equation}
where
$$ f(q^{\frac 1 2 }):=q^{-1}\frac{1}{2^{47}}
\left(\prod_{n\geq 0}(1+q^{n+\frac 1 2})+\prod_{n\geq 0}
(1-q^{n+\frac 1 2})\right)^{47}\prod_{n\geq 2}\frac{1}{(1-q^n)}.$$

Clearly, both $ch_qU$ and $f(q^{\frac 1 2})$ are convergent for
$0<|q|<1$, when $q$ is regarded as a complex number. So, we can
and do treat both $ch_qU$ and $ f(q^{\frac 1 2 })$ as functions
for $0<q<1$ and the inequality (\ref{ine}) still holds as
functions.

We have already proved in Lemma \ref{lj} that the graded dimension
of $V$ is $J(q)$ which of course also converges for $0<|q|<1.$ In
the following we will take $q$ to be a real number in the domain
$(0,1).$ Since $U$ is a subspace of $V,$ we have
$$\frac{ch_q U}{J(q)}\leq 1.$$

Let $L$ be the Niemeier lattice of type $D_{24}.$ Then the lattice
vertex operator algebra $V_L$ is a module for the affine Lie
algebra $D_{24}^{(1)}.$ Denote the irreducible highest weight
module for $D_{24}^{(1)}$ of level $k$ by $L_k(\lambda)$ where
$\lambda$ is a dominant weight of the finite dimensional Lie
algebra of type $D_{24}.$ Let $\lambda_i$ be the fundamental
weights of Lie algebra of type $D_{24}$ for $i=1,...,24$ so that
$\lambda_{23}$ and $\lambda_{24}$ are the half spin weights. (We
are using the labelling of simple roots given in \cite{H}.) Then
as a module for   $D_{24}^{(1)}$ $V_L$ is a direct sum
$$V_L=L_1(0)\oplus L_1(\lambda_{23})$$
following from the structure of lattice $L.$ It is well-known that
$$ch_qV_L=\frac{\theta_L(q)}{\eta(q)^{24}}=J(q)+2\times (24)^2-24$$
where $2\times (24)^2-24=1128$ is the dimension of the Lie algebra
of type $D_{24},$
$$\theta(q)=\sum_{\alpha\in L}q^{(\alpha,\a)/2}$$
is the theta function of the lattice $L$ and
$$\eta(q)=q^{1/24}\prod_{n\geq 1}(1-q^n).$$ So we have
$$J(q)<ch_qV_L$$
as a function in $q\in (0,1).$

On the other hand, using the Boson-Fermion correspondence given in
\cite{F},  we see that the characters of the fermion realizations
of $L_1(0)$ and $L_1(\lambda_{23})$ satisfy the following
relations
$$ch_qL_1(0)\leq ch_qL_1(0)+ ch_qL_1(\lambda_1)=q^{-1}\prod_{n\geq 0}(1+q^{n+\frac 1 2})^{48}$$
$$ch_qL_1(\lambda_{23})=q^{-1}\prod_{n>0}(1+q^n)^{48}< 2q^{-1}\prod_{n\geq 0}(1+q^{n+\frac 1 2})^{48}$$
As a result we have
$$J(q)\leq ch_q V_L\leq 3q^{-1}\prod_{n\geq 0}(1+q^{n+\frac 1 2})^{48}.$$

Note that
$$ f(q^{\frac 1 2 })\geq q^{-1}\frac{1}{2^{47}}\prod_{n\geq 0}(1+q^{n+\frac 1 2})^{47}\prod_{n\geq 2}\frac{1}{(1-q^n)}$$
So finally we have
\begin{equation}\label{fi}
\frac{ch_qU}{ch_qV}\geq \frac{1}{2^{47}3}\prod_{n\geq
0}\frac{1}{(1+q^{n+\frac 1 2})} \prod_{n\geq 2}\frac{1}{(1-q^n)}.
\end{equation}
Clearly, the right hand side of (\ref{fi}) goes to infinity as $q$
goes to 1. This is a contradiction to $\frac{ch_qU}{ch_qV}\leq 1.$
 \qed

\begin{rem} From the proof Theorem \ref{t4.2} we see that we in fact prove
a stronger result: If $\{u_1,...,u_{48}\}$ are 48 mutually
commutative Virasoro elements of central charge $\frac 1 2$ then
 $\{u_1,...,u_{48}\}$ is a VF.
% So this result has nothing to do with
%the particular choice of $\{\omega_1,...,\omega_{48}\}.$  {*** I do not understand ***}
\end{rem}

\begin{rem}\label{ra} In the proof of Theorem \ref{t4.2} we only use the fact
that $ch_q V=J(q).$ In fact, the proof goes through if we assume
that $\dim V_n\leq V^{\natural}_n$ for $n\geq 3.$ So Theorem
\ref{t4.2} holds with the assumptions given in Theorem \ref{mt2}.
\end{rem}

{\bf Proof of Theorem \ref{mt1}:}  By Theorem \ref{t4.2}, $V$ is
an FVOA with VF $F:=\{\omega_1,...,\omega_{48}\}.$ Let $U$ be the
vertex operator subalgebra generated by $V_2.$ Then $U$ is also a
FVOA with the same VF. Since $F$ is a VF in both $U$ and $V$, we
use a subscript $U$ to indicate dependence of the associated
binary codes on $U$. We have that $\sC$ is a subcode of $C_U(F)$
and $\sD$ is a subcode of $D_U(F).$ Since $D_U(F)\subset
C_U(F)^{\perp},$ and $\sD=\sC^{\perp}$ we immediately see that
$\sC=C_U(F)$ and $\sD=D_U(F).$

Note that $\sC$ is a subgroup of $C(V)$ and $\sD$ is a subgroup of
$D(F).$ Since $V$ is holomorphic by Lemma \ref{lj},
$C(F)=D(F)^{\perp}$ (see Theorem \ref{DGH}). This implies that
$C_U(F)=C(F)=\sC$ and $D_U(F)=D(F)=\sD.$

%\bsays{Earlier, elements of D were Greek letters, e.g in formulas like $M(\beta, \gamma)$, so I changed the $d$ to $\delta$ for the next few  paragraphs.  }

Now  $M_{\sC}=U^{0}$ is a vertex operator subalgebra of $V.$ Then by
Theorem \ref{DGH}, $V$ is a direct sum of inequivalent irreducible
$M_{\sC}$-modules. By Theorem \ref{indeed}, for each $\delta\in \sD$
there exists a unique $\gamma_\delta\in \Z_2^{48}/\sC$ such that
$M(\delta,\gamma_\delta)$ is isomorphic to  a submodule of $V.$ 
Then
$$V \cong \bigoplus_{\delta\in D}M(\delta,\gamma_\delta)$$
as $M_{\sC}$-module. Similarly, $V^{\natural}$ has a decomposition
$$V^{\natural} \cong \bigoplus_{\delta\in D}M(\delta,\be_\delta)$$
where $\be_\delta\in \Z_2^{48}/\sC$.    In the case that
 the lowest
weight of $M(\delta,\be_\delta)$ is 0 or 2, we have $\be_\delta=\gamma_\delta.$  Since
every module for $M_{\sC}$ is a simple current by Theorem
\ref{fusion}, by the uniqueness of simple current extension
theorem in \cite{DM2}, it is sufficient to show that
$M(\delta,\gamma_\delta)$ and $M(\delta,\be_\delta)$ are isomorphic $M_{\sC}$-modules.

For $\delta\in \sD$ we denote the lowest weight of $M(\delta,\be_\delta)$ by
$w(\delta).$ Set
$$X=\{(\delta,\be_\delta)|\, \delta\in \sD, w(\delta)=0, 2\}.$$
 Since $V^{\natural}$ is
generated by $V^{\natural}_2$ (see \cite{FLM2}),  the group
$G:=\{(\delta,\be_\delta)|\delta\in \sD\}$ is a subgroup of $D\times
\Z_2^{48}/\sC$ generated by $X.$ So, the group $H:=
\{(\delta,\gamma_\delta)|\delta\in \sD\}$ is a subgroup of $D\times
\Z_2^{48}/\sC$ and contains $G$ as a subgroup. As a result,
$G=H.$ By Theorem \ref{indeed}, $M(\delta,\gamma_\delta)$ and $M(\delta,\be_\delta)$
are indeed isomorphic $M_{\sC}$-modules. \qed

\bigskip

{\bf Proof of Theorem \ref{mt2}:} In this case, the conclusions of Lemmas \ref{l4.1},
\ref{lb}, \ref{ladd} and Theorem \ref{t4.2} still hold (see Remark
\ref{ra}).

Let $U$ be as in the proof of Theorem \ref{mt1}. Since $U$ is
generated by the Griess algebra, and $\sD\subset \sC,$  $C(U)=\sC$
and $M_{\sC}$ is a subalgebra of $U.$ From the proof of Theorem
\ref{mt1} we see that  
$$U\cong \bigoplus_{\delta\in D}M(\delta,\gamma_\delta)$$
as $M_{\sC}$-modules. The same argument used in the proof of
Theorem \ref{mt1} shows that $U$ and $V^{\natural}$ are
isomorphic. So we have
$$J(q)=ch_q V^{\natural}\geq ch_q V\geq ch_q U=J(q).$$
As a result, $U=V.$ This completes the proof. \qed

\bigskip

We give an application of Theorem \ref{mt2}. Let $U$ be the $\Z_3$
orbifold construction given in \cite{DM1}. It has been expected
for a long time that $U$ and $V^{\natural}$ are isomorphic vertex
operator algebras. The isomorphism follows from Theorem  \ref{mt2}
easily now.

\begin{cor} $V^{\natural}$ and $U$ are isomorphic.
\end{cor}

\pf $U$ satisfies the conditions in Theorem \ref{mt2}. In
particular, $\ch_q U=J(q).$ \qed

\end{document}